\documentclass[11pt,a4paper]{article}
\usepackage{graphicx}
\usepackage{subfigure}

\begin{document}

\title{\huge{Use of Triangular Elements for Nearly Exact BEM Solutions}}
\author{Supratik Mukhopadhyay, Nayana Majumdar}
\date{INO Section, Saha Institute of Nuclear Physics\\
      1/AF, Sector 1, Bidhannagar, Kolkata 700064, WB, India\\
      supratik.mukhopadhyay@saha.ac.in, nayana.majumdar@saha.ac.in}

\maketitle

\begin{abstract}
A library of C functions yielding exact solutions of potential and flux
influences due to uniform surface distribution of singularities on flat
triangular and rectangular elements has been developed. This library, ISLES, has
been used to develop the neBEM solver that is both precise and fast in solving
a wide range of problems of scientific and technological interest. Here we
present the exact expressions proposed for computing the influence of uniform
singularity distributions on triangular elements and illustrate their
accuracy. We also present a study concerning the time taken to evaluate these
long and complicated expressions \textit{vis a vis} that spent in carrying out
simple quadratures. Finally, we solve a classic benchmark problem in
electrostatics, namely, estimation of the capacitance of a unit square plate
raised to unit volt. For this problem, we present the estimated values of
capacitance and compare them successfully with some of the most accurate results
available in the literature. In addition, we present the variation of the
charge density close to the corner of the plate for various degrees of
discretization. The variations are found to be smooth and converging. This is
in clear contrast to the criticism commonly leveled against usual BEM solvers.
\end{abstract}


\section{Introduction}
One of the elegant methods for solving the Laplace / Poisson equations
(normally an integral expression of the inverse square law) is to set up the
Boundary Integral Equations (BIE) which lead to the moderately popular Boundary
Element Method (BEM). In the forward version of the BEM, surfaces of a given
geometry are replaced by a distribution of point singularities such as source /
dipole of unknown strengths. The strengths of these singularities are obtained 
through the satisfaction of a given set of boundary conditions that can be
Dirichlet, Neumann or of the Robin type. The numerical implementation requires
considerable care \cite{Nagarajan93} because it involves evaluation of singular
(weak, strong and hyper) integrals. Some of the notable two-dimensional and
three-dimensional approaches are \cite{Nagarajan93} and
\cite{Cruse69,Kutt75,Lachat76,Srivastava92,Carini2002} and the references in
these papers. Despite a large body of literature, closed form analytic
expressions for computing the effects of distributed singularities are rare
\cite{Newman86,Goto92} and complicated to implement. Thus, for solving difficult
but realistic problems involving, for example, sharp edges and corners or thin
elements, introduction of complicated mathematics and special formulations
becomes a necessity \cite{Chyuan04, Bao04}. These drawbacks are some of the
major
reasons behind the relative unpopularity of the BEM despite its significant
advantages over domain approaches such as the finite-difference and
finite-element methods (FDM and FEM) while solving non-dissipative problems
\cite{IEEE2006}. It is well-understood that most of the difficulties in the
available BEM solvers stem from the assumption of nodal concentration of
singularities which leads to various mathematical difficulties and to the
infamous numerical boundary layers \cite{Chyuan04,Sladek91}. The Inverse Square
Law Exact Solutions (ISLES) library, in contrast, is capable of truly
modeling the effect of distributed singularities precisely and, thus, its
application is not limited by the proximity of other singular surfaces or their
curvature or their size and aspect ratio. The library consists of
exact solutions for both potential and flux due to uniform distribution of
singularity on flat rectangular and triangular elements. While the rectangular
element can be of any arbitrary size \cite{EABE2006,NIMA2006}, the triangular
element can
be a right angled triangle of arbitrary size \cite{EMTM2NTriElem2007}. Since
any real geometry can be
represented through elements of the above two types (or by the triangular type
alone), this library can help in developing solvers capable of solving
three-dimensional potential problems for any geometry. It may be noted here that
any non-right-angled triangle can be easily decomposed in to two right-angled
triangles. Thus, the right-angled triangles considered here, in fact, can take
care of any three-dimensional geometry. Several difficulties were faced in
developing the library which arose due to the various terms of the integrals
and also from the approximate nature of computation in digital computers. In
this paper, we have discussed these difficulties, solutions adopted at present
and possible ways of future improvement.

The classic benchmark problem of estimating the capacitance of a unit square
plate raised to unit volt has been addressed using a solver based on ISLES,
namely, the nearly exact BEM (neBEM) solver.
Results obtained using neBEM have been compared with
other precise results available in the literature. The comparison clearly
indicates the excellent precision and efficiency achievable using ISLES and
neBEM. In addition, we have also presented the variation of charge density
close to the corner of the square plate. Usually, using BEM, it is
difficult to obtain physically consistent results close to these geometric
singularities. Wild variations in the magnitude of the charge density has been
observed with the change in the degree of discretization, the reason once
again being associated with the nodal model of singularities \cite{Wintle04}.
In contrast, using neBEM, we have obtained very smooth variation close to the
corner. Moreover, the magnitudes of the charge density have been
found to be consistently converging to physically realistic values. These
results clearly indicate that since the foundation expressions of the solver
are exact, it is possible to find the potential and flux accurately in
the complete physical domain, including the critical near-field domain using
neBEM. In addition, since singularities are no longer assumed to be nodal and
we have the exact expressions for potential and flux throughout the physical
domain, the boundary conditions no longer need to be satisfied at
\textit{special} points such as the centroid of an element. Although
consequences of this considerable advantage is still under study, it is
expected that this feature will allow neBEM to yield accurate estimates for
problems involving corners and edges that are very important in a large number
of scientific and technological studies.

It should be noted here that the exact expressions for triangular elements
consist of a
significantly larger number of mathematical operations than those for
rectangular elements. Thus, for the solver, it is more economical if we use
a mixed mesh of rectangular and triangular elements using rectangular elements
as much as possible. However, in the present work, we have intentionally
concentrated on the performance of the triangular elements and results shown
here are those obtained using only triangular elements.
\section{Exact Solutions}
The expressions for potential and flux at a point $(X, Y, Z)$ in free space due
to uniform source distributed on a rectangular flat surface having corners
situated at $(x_1, z_1)$ and $(x_2, z_2)$ has been presented, validated and used
in \cite{EABE2006, NIMA2006} and, thus, is not being repeated here.

\begin{figure}[hbt]
\begin{center}
\includegraphics[height=3in,width=3in]{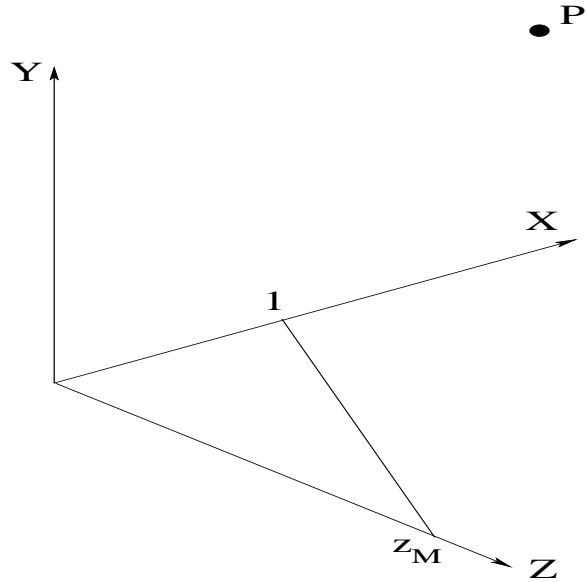}
\caption{\label{fig:GeomTriElem} Right-angled triangular element with x-length 1
and an arbitrary z-length, $z_M$; P is the point where the influence (potential
and flux) is being computed}
\end{center}
\end{figure}

Here, we present the exact expressions necessary to compute the potential and
flux due to a right-angled triangular element of arbitrary size, as shown in
Fig.\ref{fig:GeomTriElem}. It may be noted here that the length in the X
direction has been normalized, while that in the Z direction has been allowed to
be of any arbitrary magnitude, $z_M$. From the figure, it is easy to see
that in order to find out the influence due to triangular element, we have
imposed another restriction, namely, the necessity that the X and Z axes
coincide with the perpendicular sides of the right-angled triangle.
Both these restrictions are trivial and can be taken care of by carrying out
suitable scaling and appropriate vector transformations. It may be noted here
that closed-form expressions for the influence of rectangular and triangular
elements having uniform singularity distributions have been previously presented
in \cite{Newman86,Goto92}. However, in these works, the expressions presented
are quite complicated and difficult to implement. In \cite{EABE2006} and in the
present work, the expressions we have presented are lengthy, but completely
straight-forward. As a result, the implementation issues of the present
expressions, in terms of the development of the ISLES library and the neBEM
solver are managed quite easily.

It is easy to show that the influence (potential) at a point $P(X,Y,Z)$ due to
uniform source distributed on a right-angled triangular element as
depicted in Fig.\ref{fig:GeomTriElem} can be represented as a multiple of
\begin{equation}
\phi(X,Y,Z) = \int_{0}^{1} \int_{0}^{z(x)}
            \frac{dx\,dz}{\sqrt{(X-x)^2 + Y^2 + (Z-z)^2}}
\label{eqn:PotTriInt}
\end{equation}
in which we have assumed that $x_1=0$, $z_1=0$, $x_2=1$ and $z_2=z_M$, as shown
in the geometry of the triangular element. The closed-form expression for the
potential has been obtained using symbolic integration \cite{MatLabBook} which
was subsequently simplified through substantial effort. It is found to be
significantly more complicated in comparison to the expression for
rectangular elements presented in \cite{EABE2006} and can be written as
\begin{eqnarray}
\label{eqn:PotTriExact}
\lefteqn{\Phi =} \nonumber \\
&& \frac{1}{2}
	\left( \right. (z_M Y^2 - X G) (LP_1 + LM_1 - LP_2 - LM_2) 
+ i \left| Y \right| (z_M X + G) (LP_1 - LM_1 - LP_2 + LM_2) \nonumber \\
&& - S_1 X ( tanh^{-1} (\frac{R_1 + i I_1}{D_{11} |Z|})
           + tanh^{-1} (\frac{R_1 - i I_1}{D_{11} |Z|}) 
         - tanh^{-1} (\frac{R_1 + i I_2}{D_{21} |Z|})
           - tanh^{-1} (\frac{R_1 - i I_2}{D_{21} |Z|}) ) \nonumber \\
&& + i S_1 |Y| ( tanh^{-1} (\frac{R_1 + i I_1}{D_{11} |Z|})
               - tanh^{-1} (\frac{R_1 - i I_1}{D_{11} |Z|}) 
             - tanh^{-1} (\frac{R_1 + i I_2}{D_{21} |Z|})
               + tanh^{-1} (\frac{R_1 - i I_2}{D_{21} |Z|}) ) \nonumber \\
&& + \frac{2 G}{\sqrt{1+{z_M}^2}}
\log \left( \right. \frac{\sqrt{1+{z_M}^2} D_{12} - E_1}
				{\sqrt{1+{z_M}^2} D_{21} - E_2} \left. \right) 
 + 2 Z \log \left( \frac{D_{21} - X + 1}{D_{11} - X} \right) \left. \right) + C
\end{eqnarray}

where, 
\[
D_{11} = \sqrt { (X-x_1)^2 + Y^2 + (Z-z_1)^2 }; \,
D_{12} = \sqrt { (X-x_1)^2 + Y^2 + (Z-z_2)^2 }
\]
\[
D_{21} = \sqrt { (X-x_2)^2 + Y^2 + (Z-z_1)^2 }; \,
I_1 = (X-x_1)\,\left| Y \right|; \, I_2 = (X-x_2)\,\left| Y \right|
\]
\[
S_1 = {\it sign} (z_1-Z);\, R_1 = Y^2 + (Z-z_1)^2
\]
\[
E_1 = (X + {z_M}^2 - z_M Z);\, E_2 = (X - 1 - z_M Z), \,
\]
\[
G = z_M (X - 1) + Z; \,
H_1 = Y^2 + G (Z - z_M);\, H_2 = Y^2 + G Z
\]
\begin{eqnarray*}
LP_1 =
	 \frac{1}{G - i z_M |Y|} 
	log(\frac{(H_1\, + G D_{12}) + i |Y| (E_1 - i z_M D_{12})} {- X + i |Y|})
\end{eqnarray*}
\begin{eqnarray*}
LM_1 = 
   \frac{1}{G + i z_M |Y|} 
	log(\frac{(H_1\, + G D_{12}) - i |Y| (E_1 - i z_M D_{12})} {- X - i |Y|})
\end{eqnarray*}
\begin{eqnarray*}
LP_2 = 
   \frac{1}{G - i z_M |Y|} 
	log(\frac{(H_2\, + G D_{21}) + i |Y| (E_2 - i z_M D_{21})} {1 - X + i |Y|})
\end{eqnarray*}
\begin{eqnarray*}
LM_2 = 
   \frac{1}{G + i z_M |Y|} 
	log(\frac{(H_2\, + G D_{21}) - i |Y| (E_2 - i z_M D_{21})} {1 - X - i |Y|})
\end{eqnarray*}
and $C$ denotes a constant of integration.

Similarly, the flux components due to the above singularity distribution can
also be represented through closed-form expressions as shown below:

\begin{eqnarray}
\label{eqn:FxTriExact}
\lefteqn{F_x = -\frac{\partial \Phi}{\partial x} =} \nonumber \\
&& \frac{1}{2}
	\left( \right. (G) (LP_1 + LM_1 - LP_2 - LM_2) 
 - i \left| Y \right| (z_M) (LP_1 - LM_1 - LP_2 + LM_2) \nonumber \\
&& + S_1 ( tanh^{-1} (\frac{R_1 + i I_1}{D_{11} |Z|})
           + tanh^{-1} (\frac{R_1 - i I_1}{D_{11} |Z|}) 
         - tanh^{-1} (\frac{R_1 + i I_2}{D_{21} |Z|})
           - tanh^{-1} (\frac{R_1 - i I_2}{D_{21} |Z|}) ) \nonumber \\
&& + \frac{2 z_M}{\sqrt{1+{z_M}^2}}
\log \left( \right. \frac{\sqrt{1+{z_M}^2} D_{12} - E_1}
				{\sqrt{1+{z_M}^2} D_{21} - E_2} \left. \right) \left. \right)
				+ C
\end{eqnarray}

\begin{eqnarray}
\label{eqn:FyTriExact}
\lefteqn{F_y = -\frac{\partial \Phi}{\partial y} =} \nonumber \\
&& \frac{-1}{2}
	\left( \right. (2 z_M Y) (LP_1 + LM_1 - LP_2 - LM_2) 
 + i \left| Y \right| (Sn(Y) G) (LP_1 - LM_1 - LP_2 + LM_2) \nonumber \\
&& + i S_1 Sn(Y) ( tanh^{-1} (\frac{R_1 + i I_1}{D_{11} |Z|})
               - tanh^{-1} (\frac{R_1 - i I_1}{D_{11} |Z|}) 
             - tanh^{-1} (\frac{R_1 + i I_2}{D_{21} |Z|})
         + tanh^{-1} (\frac{R_1 - i I_2}{D_{21} |Z|}) ) \left. \right) 
			+ C \nonumber \\
&&
\end{eqnarray}

and,

\begin{eqnarray}
\label{eqn:FzTriExact}
F_z = -\frac{\partial \Phi}{\partial z} =
\left( \right . \frac{1}{\sqrt{1+{z_M}^2}}
\log \left( \right. \frac{\sqrt{1+{z_M}^2} D_{21} - E_2}
				{\sqrt{1+{z_M}^2} D_{12} - E_1} \left. \right) 
 + \log \left( \frac{D_{11} - X}{D_{21} - X + 1} \right) \left. \right) + C
\end{eqnarray}
where $Sn(Y)$ implies the sign of the Y-coordinate and $C$ indicates constants
of integrations. It is to be noted that the constants of different integrations
are not the same.
These expression are expected to be useful in the mathematical modeling of
physical processes governed by the inverse square laws. Being exact and valid
throughout the physical domain, they can be used to formulate versatile solvers
to solve multi-scale multi-physics problems governed by the Laplace / Poisson
equations involving Dirichlet, Neumann or Robin boundary conditions.

\section{Development of the ISLES library}
Due to the tremendous popularity of the C language
we have written the codes in the C programming language. However, it should be
quite simple to translate the library to other popular languages such as
FORTRAN or C++, since no special feature of the C language has been used to
develop the codes.

\subsection{Validation of the exact expressions}
\begin{figure}[hbt]
\begin{center}
\includegraphics[height=4in,width=4in]{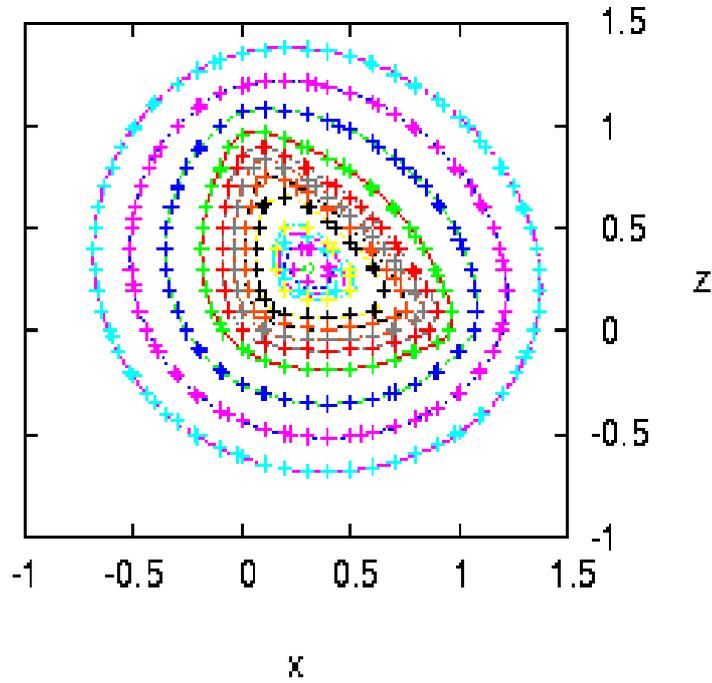}
\caption{\label{fig:CntrTriElem} Potential contours on a triangular element computed using exact expressions and by numerical quadrature}
\end{center}
\end{figure}

\begin{figure}[hbt]
\begin{center}
\includegraphics[height=3in,width=5in]{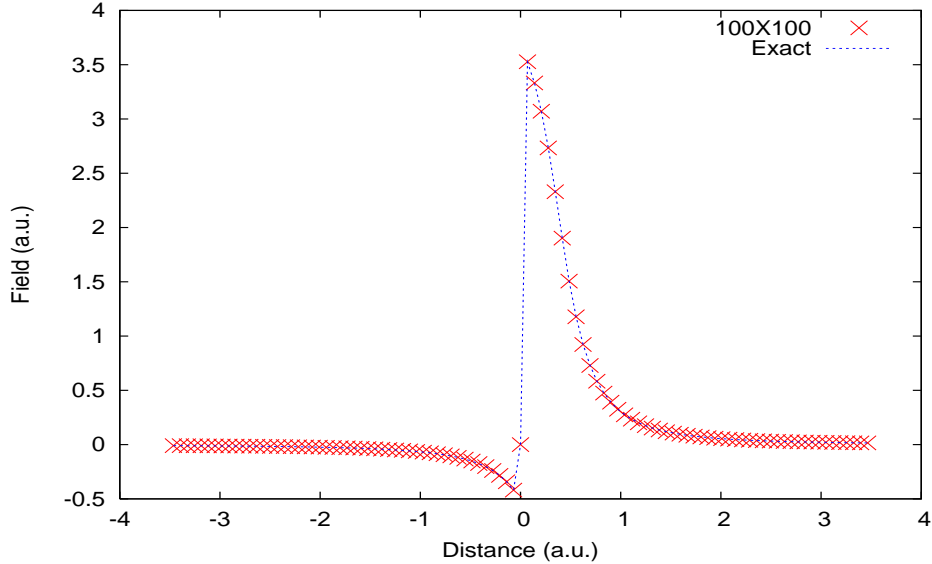}
\caption{\label{fig:TriDiag1FY} Comparison of flux (in the Y direction) as
computed by ISLES and numerical quadrature along a diagonal line}
\end{center}
\end{figure}

\begin{figure}[hbt]
\centering
\subfigure[Potential surface]{\label{fig:SurfTriElem}\includegraphics[width=0.45\textwidth]{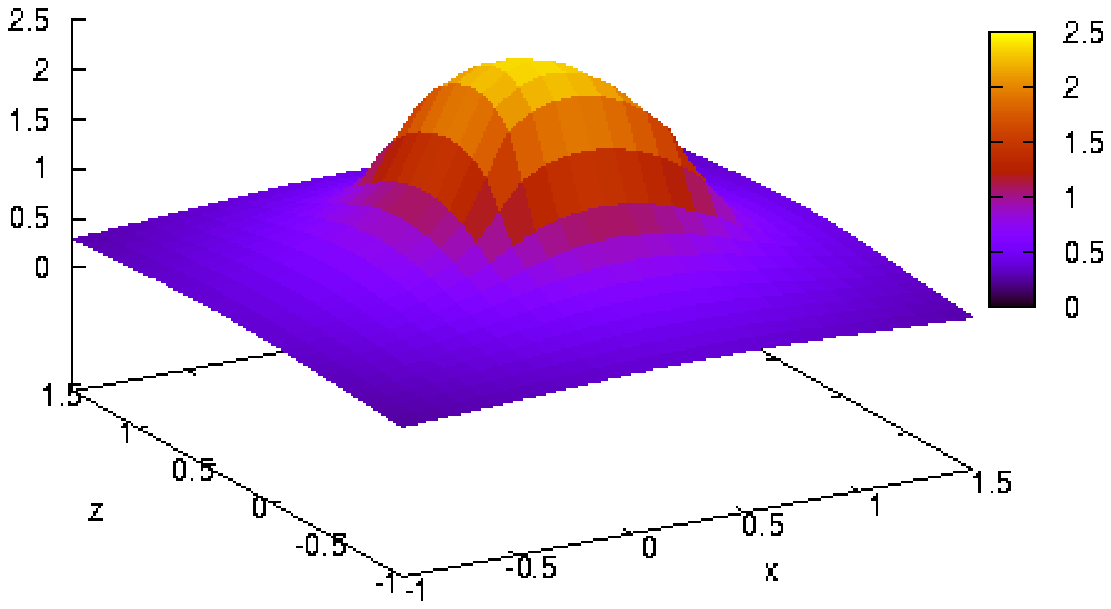}}
\subfigure[Flux surface]{\label{fig:TriFY_XY0}\includegraphics[width=0.45\textwidth]{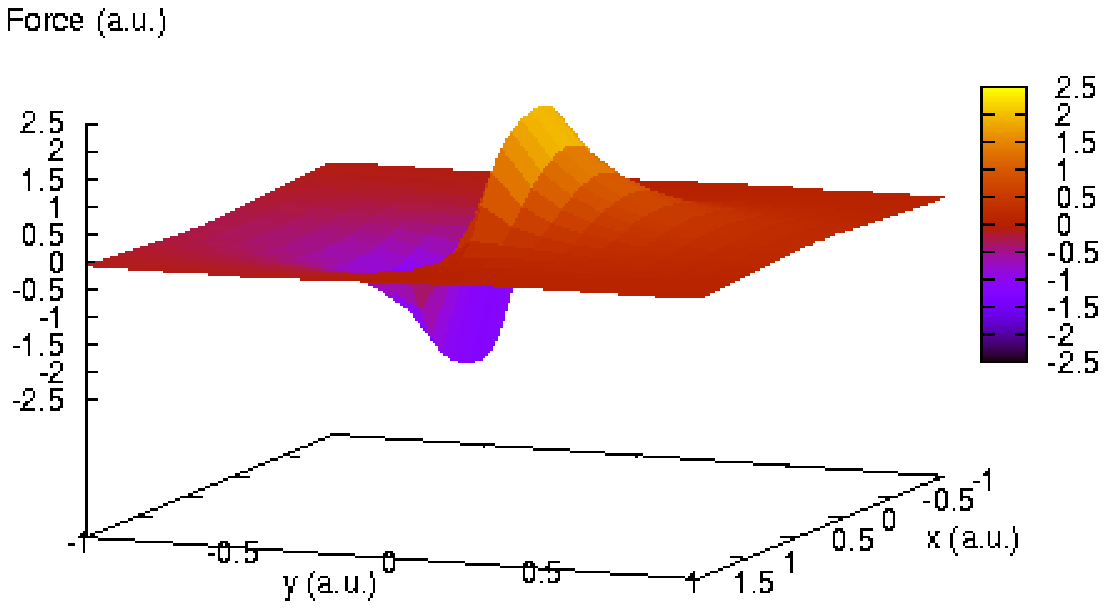}}
\caption{(a) Potential surface due to a triangular source distribution on the
element plane,
(b) Flux (in the Y direction) surface due to a triangular source distribution
on the XY plane at Z=0}
\label{fig:SurfacePlots}
\end{figure}

The expressions for the rectangular element have been validated in detail in
\cite{EABE2006}. Here, we present the results for triangular elements in fair
detail. In Fig.\ref{fig:CntrTriElem}, we have presented a comparison of
potentials evaluated for a unit triangular element by using the exact
expressions, as well as by using numerical quadrature of high accuracy. The two
results are found to compare very well throughout. Please note that contours
have been obtained on the plane of the element, and thus, represents a rather
critical situation. Similarly, Fig.\ref{fig:TriDiag1FY}
shows a comparison between the results obtained using closed-form expressions
for flux and those obtained using numerical quadrature. The flux considered here
is in the $Y$ direction and is along a line beginning from $(-2,-2,-2)$ and
ending at $(2,2,2)$. The comparison shows the commendable accuracy expected from
closed form expressions.
In Fig.\ref{fig:SurfTriElem} and \ref{fig:TriFY_XY0}, the surface plots of
potential on the element plane ($XZ$ plane) and $Y$-flux on the $XY$ plane
have been presented from which the expected significant increase
in potential and sharp change in the flux value on the element is observed.
Thus, by using a small
fraction of computational resources in comparison to those consumed in numerical
quadratures, ISLES can compute the exact value of potential and flux
for singularities distributed on triangular elements.

\subsection{Near-field performance}
\begin{figure}[hbt]
\begin{center}
\includegraphics[height=3in,width=5in]{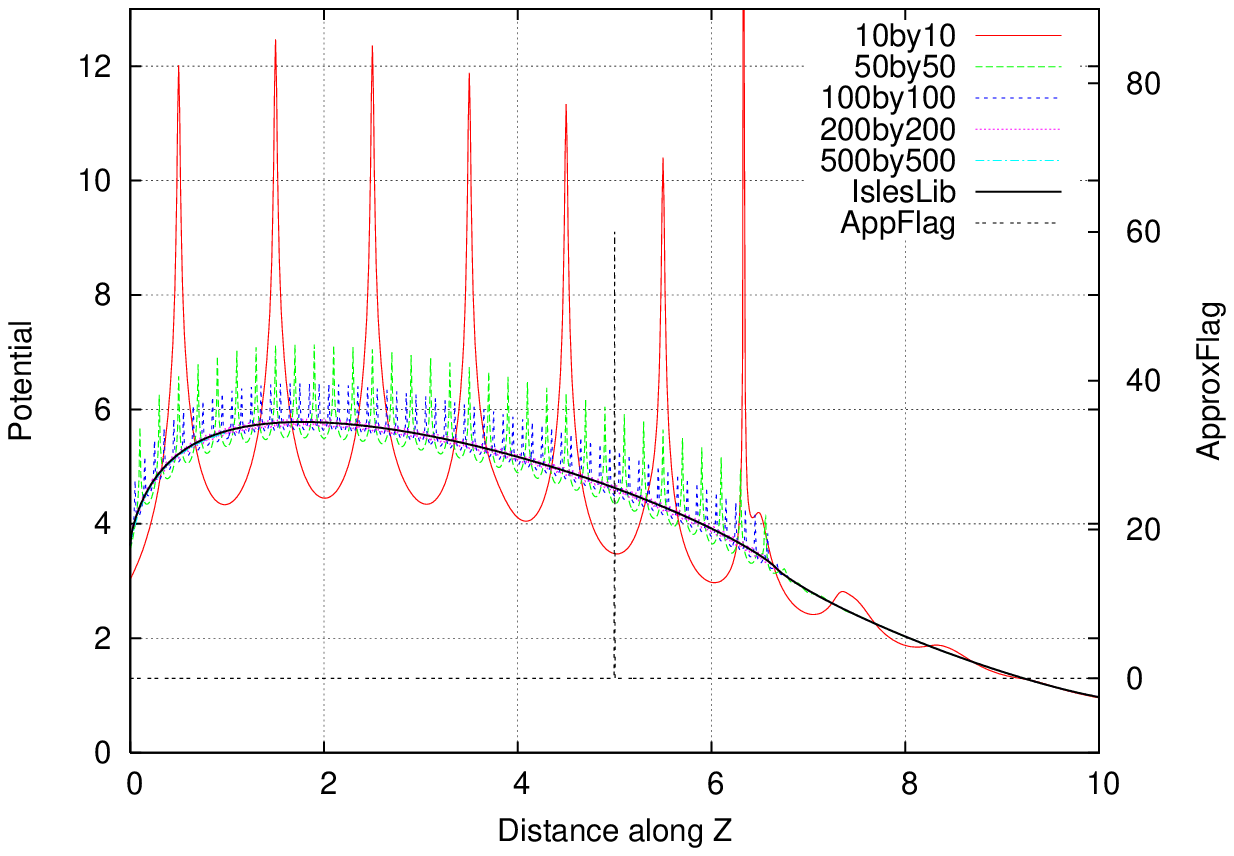}
\caption{\label{fig:TCentroidZVsPot} Variation of potential along a centroidal
line on the XZ plane parallel to the Z axis for a triangular element: comparison
among values obtained using the exact expressions and numerical quadratures}
\end{center}
\end{figure}

\begin{figure}[hbt]
\begin{center}
\includegraphics[height=3in,width=5in]{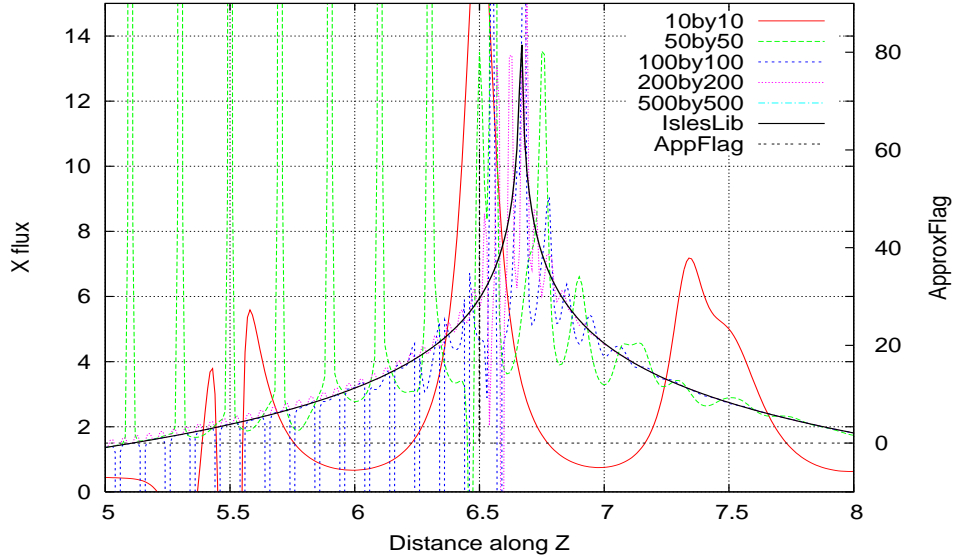}
\caption{\label{fig:TCentroidZVsFx} Variation of flux in the X direction along a
line on the XZ plane parallel to the Z axis for a triangular element: comparison
among values obtained using the exact expressions and numerical quadratures}
\end{center}
\end{figure}

\begin{figure}[hbt]
\begin{center}
\includegraphics[height=3in,width=5in]{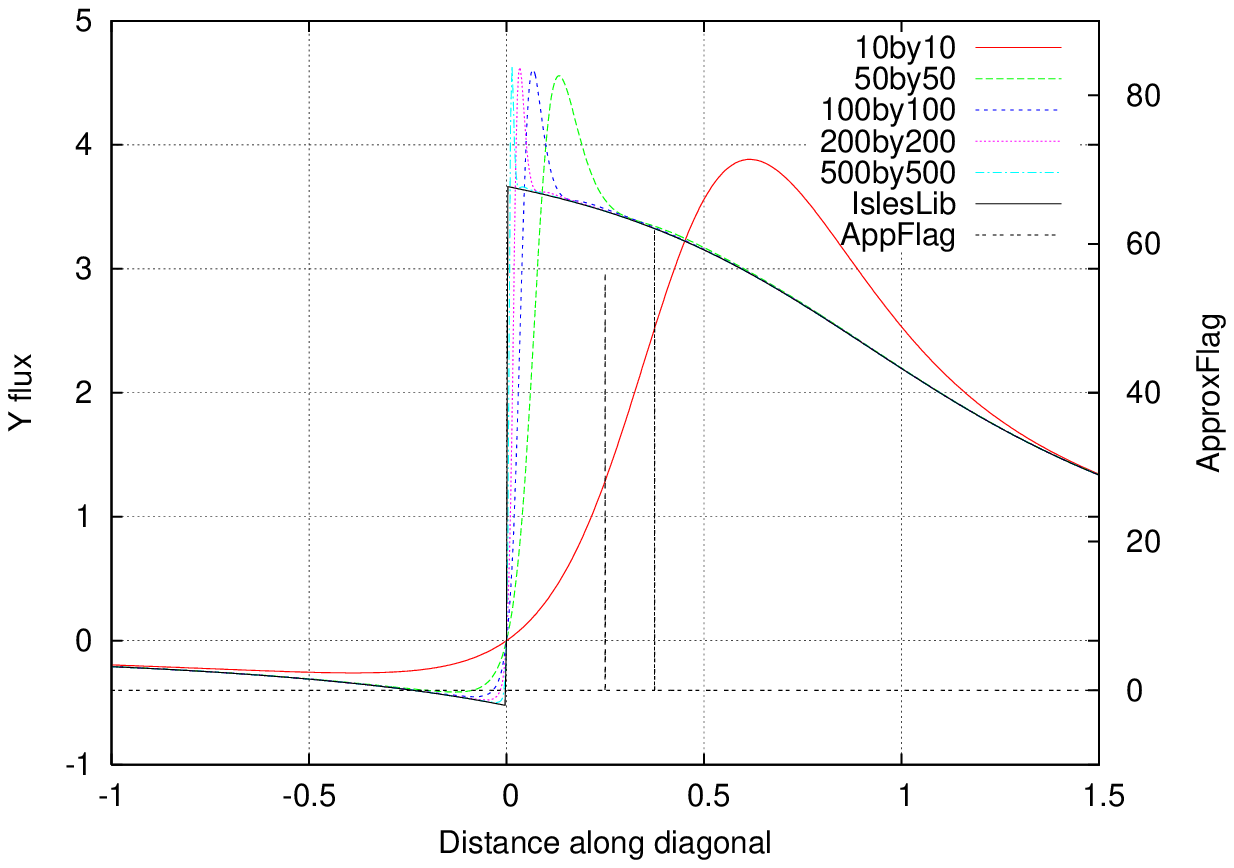}
\caption{\label{fig:TDiag1VsFy} Comparison of flux (in the Y direction) along a
diagonal line piercing the triangular element at the centroid: comparison
among values obtained using the exact expressions and numerical quadratures}
\end{center}
\end{figure}

In order to emphasize the accuracy of ISLES, we have considered the following
severe situations in the near-field region in which it is observed that the
quadratures can match the
accuracy of ISLES only when a high degree of discretization is used. Please
note that in these cases, the value of $z_M$ has been considered to be 10. In
Fig.\ref{fig:TCentroidZVsPot} we have presented the variation of potential
along a line on the element surface running parallel to the Z-axis of the
triangular element
(see Fig.\ref{fig:GeomTriElem}) and going through the centroid of the element.
It is observed that results obtained using even a $100 \times 100$ quadrature is
quite unacceptable. In fact, by zooming on to the image, it can be found that
only the maximum discretization yields results that match closely to the
exact solution. It may be noted here that the potential is a relatively easier
property to compute. The difficulty of achieving accurate flux estimates is
illustrated in the two following figures. The variation of flux in the
$X$-direction along the same line as used in Fig.\ref{fig:TCentroidZVsPot} has
been presented in Fig.\ref{fig:TCentroidZVsFx}. Similarly, variation of $Y$-flux
along a diagonal line (beginning at (-10,-10,-10) and ending at (10,10,10)
and piercing the element at the centroid) has been presented in
Fig.\ref{fig:TDiag1VsFy}. From these figures we see that the flux values
obtained using the quadrature are always inaccurate even if the discretization
is as high as $100 \times 100$.
We also observe that the estimates are locally inaccurate
despite the use of very high amount of discretization ($200 \times 200$ or
$500 \time 500$). Specifically, in the latter figure, even the highest
discretization can not match the exact values at the peak, while in the former
only the highest one can correctly emulate the sharp change in the flux value.
It is also heartening to note that the values from the quadrature using higher
amount of discretization consistently converge towards the ISLES values.

\subsection{Far field performance}
\label{sec:FarField}
It is expected that beyond a certain distance, the effect of the singularity
distribution can be considered to be the same as that of a centroidally
concentrated singularity or a simple quadrature. The optimized amount of
discretization to be used for the quadrature can be determined from a study of
the speed of execution of each of the functions in the library and has been
presented separately in a following sub-section. If we plan to replace the
exact expressions by quadratures (in order to reduce the computational expenses,
presumably) beyond a certain given distance, the quadrature should
necessarily be efficient enough to justify the replacement. While standard but
more elaborate algorithms similar to the fast multipole method (FMM)
\cite{Greengard87} along with the GMRES \cite{Saad86} matrix solver can lead to
further of computational efficiency, the simple approach as outlined
above can help in reducing a fair amount of computational effort. In the
following, we present the results of numerical experiments that help us in
determining the far-field performance of the exact expressions and quadratures
of various degrees that, in turn, help us in choosing the more efficient
approach for a desired level of accuracy.

\begin{figure}[hbt]
\begin{center}
\includegraphics[height=3in,width=5in]{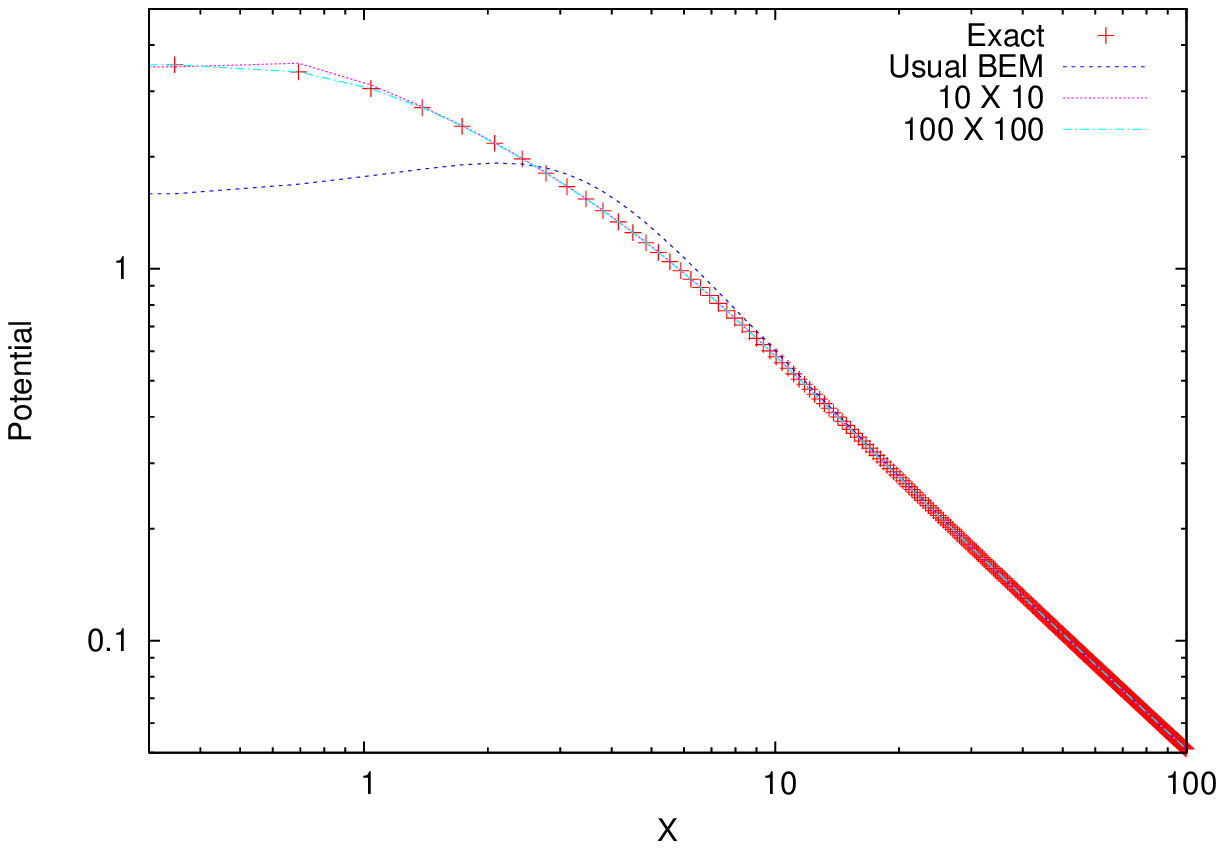}
\caption{\label{fig:DiagPFFarField} Potential along a diagonal through the
triangular element computed using exact, $100 \times 100$,  $10 \times 10$ and
usual BEM approach}
\end{center}
\end{figure}

\begin{figure}[hbt]
\begin{center}
\includegraphics[height=3in,width=5in]{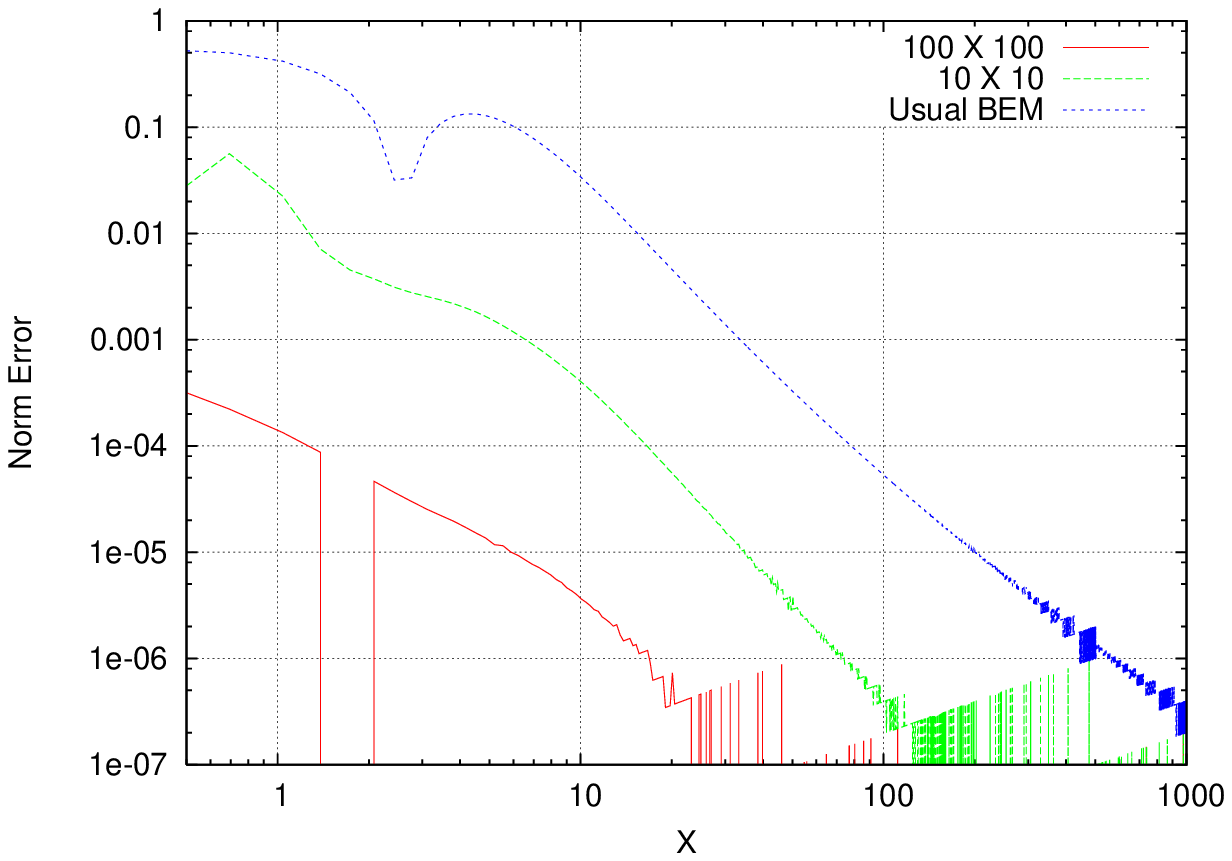}
\caption{\label{fig:ErrorPot} Error along a diagonal through the triangular
element computed using $100 \times 100$, $10 \times 10$ and usual BEM approach}
\end{center}
\end{figure}

In Fig.\ref{fig:DiagPFFarField} we have presented potential values obtained
using the exact approach, $100 \times 100$, $10 \times 10$ and no
discretization, i.e., the
usual BEM approximation while using the zeroth order piecewise uniform charge
density assumption. The potentials are computed along a diagonal line running
from (-1000, -1000, -1000) to (1000, 1000, 1000) which pierces a triangular
element of $z_M = 10$. It can be seen that results obtained using the usual BEM
approach yields inaccurate results as we move closer than distances of 10 units,
while the $10 \times 10$ discretization yields acceptable results up to a
distance
of 1.0 unit. In order to visualize the errors incurred due to the use of
quadratures, we have plotted Fig.\ref{fig:ErrorPot} where the errors
incurred (normalized with respect to the exact value) have been plotted. From
this figure we can conclude that for the
given diagonal line, the error due to the usual BEM approximation falls below
1\% if the distance is larger than 20 units while for the simple $10 \times 10$
discretization, it is 2 units. It may be mentioned here that along the axes the
error turns out to be significantly more \cite{EABE2006} and the limits
need to
be effectively doubled to achieve the accuracy for all cases possible. Thus,
for achieving 1\% accuracy, the usual BEM is satisfactory only if the distance
of the influenced point is five times the longer side of an element. Please note
here that the error drops to 1 out of $10^6$ as the distance becomes fifty times
the longer side. Besides proving that the exact expressions work equally
well in the near-field as well as the far-field, this fact justifies the usual
BEM approach for much of the computational domain leading to substantial savings
in computational expenses.

The accuracy of the exact expressions used in the ISLES library is confirmed
from the above comparisons. However, there are several other important issues
related to the development of the library that are discussed below briefly.

\subsection{Evaluation of the component functions}
Many of the irrational and transcendental functions have domains and ranges
in which they are defined. Moreover, they are often multiply defined in the
complex domain; for example, there are an infinite number of complex values for
the logarithm function. In such cases, only one principal value must be returned
by the function. In general, such values cannot be chosen so as to make the
range continuous and thus, lines in the domain called branch cuts need to be
defined, which in turn define the discontinuities in the range. While evaluating
expressions such as the ones displayed in eqns.(\ref{eqn:PotTriExact}
- \ref{eqn:FzTriExact}) a number of such problems are
expected to occur. However, when the expressions are analyzed at critical
locations such as the corners and edges of the element, it is observed that
the terms likely to create difficulties while evaluating potentials are either
cancelled out or are themselves multiplied by zero. As a result, at these
locations of likely geometric and mathematical singularities, the solution
behaves nicely. However, the same is not true for the expressions related to the
flux components. For these, we have to deal with branch-cut problems in relation
to $tanh^{-1}$ and problems related to the evaluation of $log(0)$. It should be
noted here that these singularities associated to the edges and corners of the
elements are of the \textit{weak} type and it is expected that exact evaluation
of these terms as well will be possible through further work.

However, difficulties of a different nature crop up in these calculations which
can be linked directly to the limitation of the computer itself, namely,
round-off errors \cite{EMTM2NTriElem2007}. These errors can
lead to severe problems while handling multi-scale problems such as those
described in \cite{EABE2006}. A completely
different approach is necessary to cope up with these difficulties, for example,
the use of extended range arithmetic \cite{Smith81}, interval arithmetic
\cite{Alefeld83} or the use of specialized
libraries such as the CORE library of the Exact Geometric Computation (EGC)
initiative \cite{EGC}. 
In the present version of ISLES, a simple approach has been implemented which
sets a lower limit to various distance values. Below this value, the distance
is considered to be zero. Plan of future improvements
in this regard has been kept at a high priority.

\subsubsection{Algorithm}
As discussed above, there are possibilities of facing problems while using the
exact expressions which may be due to
the functions being evaluated or due to round-off errors leading to erroneous
results. Moreover, despite providing many checks during the computation
there is finite possibility of ending up with a wrong value of a property
indicated by its being \textit{Nan} or \textit{inf} or potential due to unit
positive singularity strength turning out to be \textit{negative}. In order
to maintain the robustness of the library, we have tried to keep checks on the
intermediate and final values during the course of the computation. When the
results are found to be unsatisfactory, unphysical, we have re-estimated the
results by using numerical quadrature and kept a track of the cause by raising
a unique approximation flag which is specific for a problem. As a result, the
steps for the calculation for a property can be written as follows:
\begin{itemize}
\item
{Get the required inputs - geometry of the element and the position where the
effect needs to be evaluated; Check whether the element size and distances
are large enough so that the results do not suffer from round-off errors.}
\item
{Check whether the location coincides with one of the \textit{special} ones,
such as corners or edges.}
\item
{Evaluate the necessary expressions in accordance with
the foregoing results. If necessary, consider each term in the expressions
separately to sort out difficulties related to singularities, branch-cuts or
round-off errors. Note that if the multiplier is zero, rest of the term does
not need evaluation.}
\item
{If direct evaluation of the expressions fail, raise a unique approximation flag
specific to this problem and term and return the value of the property by using
numerical quadrature.}
\item
{Compute all the terms and find the final value, Check whether the final value
is a number and physically meaningful. If not satisfied, recompute the result
using numerical quadrature and raise the relevant approximation flag.}
\end{itemize}

\subsection{Speed of execution}
The time taken to compute the potential and flux is an important parameter
related to the overall computational efficiency of the codes. This is true
despite the fact that, in a typical simulation, the time taken to solve the
system of algebraic equation is far greater than the time taken to build the
influence coefficient matrix and post-processing. Moreover, the amount of time
taken to solve the system of equations tend to increase at a greater rate than
the time taken to complete the other two.
It should be mentioned here that the time taken in each of these steps can
vary to a significant amount depending on the algorithm of the solver. In the
present case, the system of equations has been solved using lower upper
decomposition using the well known Crout's partial pivoting. Although this
method is known to be very rugged and accurate, it is not efficient as far as
number of arithmetic operations, and thus, time is concerned. It is also
possible to reduce the time taken to pre-process (generation of mesh and
creation of influence matrices), solve the system of algebraic equations and
that for post-process (computation of potential and flux at the required
locations) can be significantly reduced by adopting faster algorithms, including
those involving parallelization. 

In order to optimize the time taken to generate the influence coefficient matrix
and that to carry out the post-processing, we carried out a small numerical
study to determine the amount of time taken to complete the various functions
being used in ISLES, especially those being used to evaluate the exact
expressions and those being used to carry out the quadratures. The results of
the study (which was carried out using the linux system command \textit{gprof})
has been presented in the following Table\ref{Table:SpeedTable}.

\begin{table*}[hbt]
\centering
\caption{\label{Table:SpeedTable}Time taken to evaluate exact expression and
various quadratures}
\begin{tabular}{| l | c | c | c | c | c |}
\hline
Method & Exact       & Usual BEM & $10 \times 10$ & $100 \times 100$ & $500 \times 500$ \\
\hline
Time   & 0.8 $\mu s$ & 25 $ns$   & 1 $\mu s$      & 200 $\mu s$      & 5 $ms$ \\
\hline
\end{tabular}
\end{table*}

Please note that the numbers presented in this table are representative and
are likely to have statistical fluctuations. However, despite the fluctuations,
it may be safely concluded that a quadrature having only $10 \times 10$
discretization is already consuming time that is comparable to that needed exact
evaluation. Thus, the exact expressions, despite their complexity, are extremely
efficient in the near-field which can be considered at least as large as 0.5
times the larger side of a triangular element (please refer to
Fig.\ref{fig:ErrorPot}). In making this statement, we have assumed that the
required accuracy for generating the influence coefficient matrix and subsequent
potential and flux calculations is 1\%. This may not be acceptable at all under
many practical circumstances, in which case the near-field would imply a larger
volume.

\subsection{Salient features of ISLES}
Development of usual BEM solvers are dependent on the two following assumptions:
\begin{itemize}
\item {While computing the influences of the singularities, the singularities
are modeled by a sum of known basis functions with constant unknown
coefficients. For example, in the constant element approach, the singularities
are assumed to be concentrated at the centroid of the element, except for
special cases such as self influence. This becomes necessary because
closed form expressions for the influences are not, in general, available for
surface elements. An approximate and computationally rather expensive way of
circumventing this limitation is to use numerical integration over each
element or to use linear or higher order basis functions.}
\item{The strengths of the singularities are solved depending upon the boundary
conditions, which, in turn, are modeled by the shape
functions. For example, in the constant element approach, it is assumed that
it is sufficient to satisfy the boundary conditions at the centroids of
the elements. In this approach, the position of the singularity and the point
where the boundary condition is satisfied for a given element usually matches
and is called the collocation point.}
\end{itemize}
The first (and possibly, the more damaging) approximation for BEM solvers can
be relaxed by using ISLES and can be restated as,
\begin{itemize} 
\item {The singularities distributed on the boundary elements are assumed to be
uniform on a particular element. The strength of the singularity may change from
element to element.}
\end{itemize}
This improvement turns out to be very significant as demonstrated in the
following section and some of our other studies involving microelectromechanical
systems (MEMS) and gas detectors for nuclear applications \cite{EABE2006,
NIMA2006}. Some of the advantages of using ISLES are itemized below:
\begin{itemize}
\item{For a given level of discretization, the estimates are more accurate,}
\item{Effective efficiency of the solver improves, as a result,}
\item{Large variation of length-scales, aspect ratios can be tackled,}
\item{Thinness of members or nearness of surfaces does not pose any problem},
\item{Curvature has no detrimental effect on the solution,}
\item{The boundary condition can be satisfied anywhere on the elements,
i.e., points other than the centroidal points can be easily used, if necessary
(for a corner problem, may be),}
\item{The same formulation, library and solver is expected to work in majority
of physical situations. As a result, the necessity for specialized formulations
of BEM can be greatly minimized.}
\end{itemize}

\section{Capacitance of a unit square plate - a classic benchmark problem}
Using the neBEM solver, we have computed the capacitance of a unit square
conducting plate raised to a unit volt. This problem is still considered to be
one of the major unsolved problems of electrostatic theory
\cite{Maxwell,Goto92,Read,Wintle04}
and no analytical solution for this problem has been obtained so far.
The capacitance value estimated by the present method has been compared with
very accurate results available in the literature (using BEM and other methods).
The results obtained using the neBEM solver is found to be among the most
accurate ones available till date as shown in Table.\ref{table:CapComp}. Please
note that we have not invoked symmetry or used extrapolation techniques to
arrive at our result presented in the table.

\begin{table*}[hbt]
\centering
\caption{\label{table:CapComp}Comparison of capacitance values}
\begin{tabular}{| l | l | c |}
\hline
Reference & Method & Capacitance (pF) / 4 $\pi \epsilon_0$ \\
\hline
\cite{Maxwell} & Surface Charge & 0.3607 \\
\hline
\cite{Reitan} & Surface Charge & 0.362 \\
\hline
\cite{Solomon} & Surface Charge & 0.367 \\
\hline
\cite{Goto92} & Refined Surface Charge & $0.3667892 \pm 1.1 \times 10^{-6}$ \\
& and Extrapolation & \\
\hline
\cite{Read} & Refined Boundary Element & $0.3667874 \pm 1 \times 10^{-7}$ \\
& and Extrapolation & \\
\hline
\cite{Mansfield} & Numerical Path Integration & $0.36684$ \\
\hline
\cite{Wintle04} & Random Walk & $0.36 \pm 0.01$ \\
\hline
This work & neBEM & 0.3660587 \\
\hline
\end{tabular}
\end{table*}

Finally, we consider the corner problem related to the electrostatics of the
above conducting plate. Problems of this nature are considered to be
challenging for any numerical tool and especially so for the BEM approach.
The inadequacy of the BEM approach, especially in solving the present problem,
has been mentioned even quite recently \cite{Wintle04}
in which it has been correctly mentioned that since the method can not extend
its mathematical model to include the edges and corners in reality, it is
unlikely that it will ever succeed in modeling the edge / corner singularities
correctly. As a result, with change in discretization, the properties near these
geometric singularities are expected to oscillate significantly leading to
erroneous results. However, as discussed above, the neBEM does extend its
singularities
distributed on the surface elements right till an edge or a corner. Moreover,
using neBEM, it is also possible to satisfy the boundary conditions (both
potential and flux) as close to the edge / corner as is required. In fact, it
should be possible to specify the potentials right at the edge / corner.

\begin{figure}[hbt]
\begin{center}
\includegraphics[height=3in,width=5in]{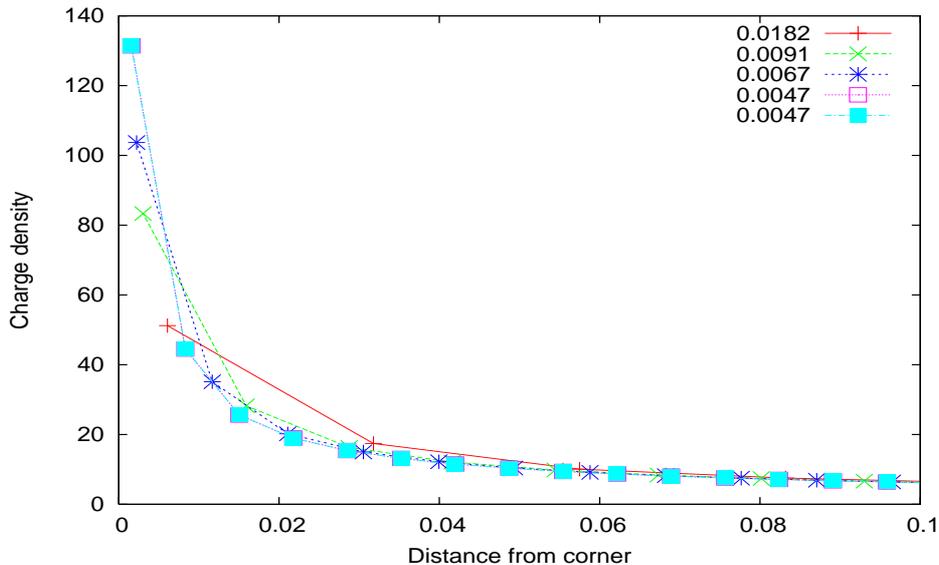}
\caption{\label{fig:CompareCornerChDen} Corner charge density estimated by 
\textit{neBEM} using various sizes of triangular elements}
\end{center}
\end{figure}

In the following study, we have presented estimates of charge density very close
the flat plate corner as obtained using neBEM. Please note that the boundary
conditions have been satisfied at the centroids of each element although we
plan to carry out detailed studies of changing the position of these points,
especially in relation to problems involving edges / corners. In
Fig.\ref{fig:CompareCornerChDen}, charge densities very close to the corner
of the flat plate estimated by neBEM using various amounts of discretization
have been presented. It can be seen that each curve follows the same general
trend, does not suffer from any oscillation and seems to be converging to a
single
curve. This is true despite the fact that there has been almost an order of
magnitude variation in the element lengths.

\begin{figure}[hbt]
\begin{center}
\includegraphics[height=3in,width=5in]{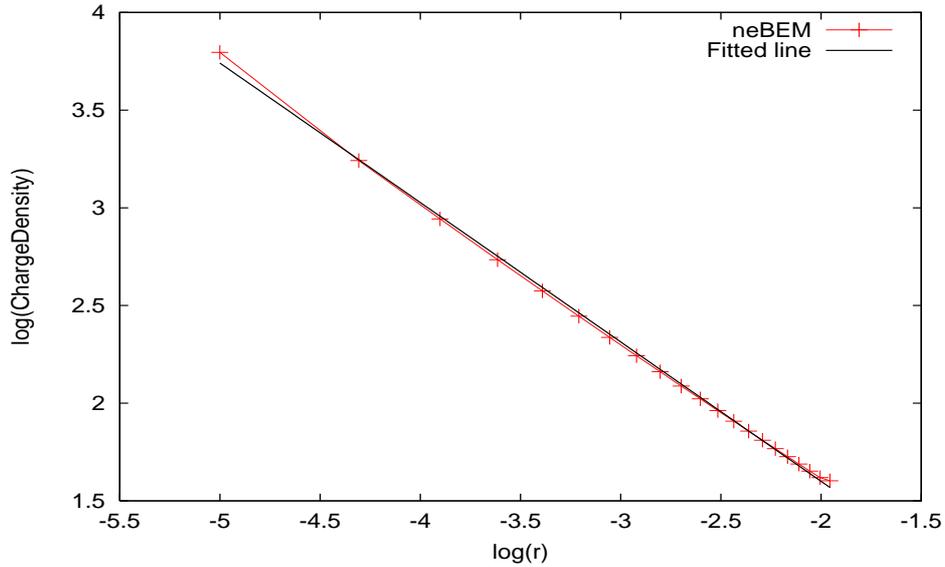}
\caption{\label{fig:FitCornerChDen} Variation of charge density with increasing
distance from the corner of the unit square plate and a least-square fitted
straight line: slope of the fitted line is 0.713567}
\end{center}
\end{figure}

Finally, in Fig.\ref{fig:FitCornerChDen}, we present a least-square fitted
straight line matching the charge density as obtained the highest discretization
in this study. It is found that the slope of the straight line is 0.713567,
which compares very well with both old and recent estimates of 0.7034
\cite{Morrison75, Hwang2005}. This is despite the fact that here we have used
a relatively coarse discretization near the corner. It should be mentioned
that none of the earlier references cited here used the BEM approach. While the
former used a singular perturbation technique, the latter
used a diffusion based Monte-Carlo method. Thus, it is extremely encouraging
to note that using the neBEM approach, we have been able to match the accuracy
of these sophisticated techniques.


\section{Conclusion}
An efficient and robust library for solving potential problems in a large
variety of science and engineering problems has been developed. Exact
closed-form expressions used to develop ISLES have been validated throughout
the physical domain (including the critical near-field region) by comparing
these results with results obtained using numerical quadrature of high
accuracy. Algorithmic aspects of this development have also been touched upon.
A classic benchmark problem of electrostatics
has been successfully simulated to very high precision.
Charge density values at critical geometric locations like corners
have been found to be numerically stable and physically acceptable.
Several advantages over usual BEM solvers and other specialized BEM solvers
have been briefly mentioned.
Work is under way to make the code more robust and efficient through
the implementation of more efficient algorithms and parallelization.

\vspace{0.25in}
\textbf{Acknowledgements}\\
We would like to thank Professor Bikas Sinha, Director, SINP and Professor
Sudeb Bhattacharya, Head, INO Section, SINP for their support and encouragement
during the course of this work.

%
%

\end{document}